\documentstyle[twoside,12pt]{article}
\input epsf
\def\bee{\begin{equation}}
\def\ee{\end{equation}}
\def\tytul#1{
        \ifodd \value{page} {\Bx\protect\hfill{#1} \vskip 0.4cm} \else {\Bx
\noindent{#1}\vskip 0.3cm}\fi
}

\begin{document}

\thispagestyle{empty}
\bigskip\bigskip
\centerline{    }
\vskip 4 cm
\centerline{\Large\bf  Some Remarks on the Distribution of twin Primes}
\bigskip\bigskip\bigskip
\centerline{\large\sl Marek Wolf}
\bigskip
\centerline{\sf Institute of Theoretical Physics, University of Wroc{\l}aw}
\centerline{\sf Pl.Maxa Borna 9, PL-50-204 Wroc{\l}aw, Poland}
\centerline{\small e-mail:  mwolf@ift.uni.wroc.pl}

\bigskip\bigskip\bigskip
\begin{abstract}

The computer data up to $2^{44}\approx 1.76\times 10^{13}$ on the gaps
between consecutive twins is presented.
The simple derivation of the heuristic formula describing
computer results contained in the recent papers by P.F.Kelly and T.Pilling
\cite{Kelly1}, \cite{Kelly2} is provided and compared with the
``experimental'' values.

\end{abstract}



\vfill
\eject

\pagestyle{myheadings}

Among the primes the subset of twin
primes is distinguished: twins are such numbers $(p, p^\prime)$ that
both $p$ and $p^\prime=p+2$ are prime. So the set of twins starts with
(3, 5), (5, 7), (11, 13), (17, 19), (29, 31) .... It is not known
whether there is
infinity of twins; the largest known today pair of twins  was found
recently by Underbakke, Carmody, Gallot (see
http://www.utm.edu/research/primes/largest.html) and  it is  the
following pair of 29603 digits numbers:
\bee
1807318575\times 2^{98305}\pm 1.
\ee

The mathematicians are using the notation $\pi_2(N)$
to denote the number of twins smaller than $N$ and the
Hardy and Littlewood conjecture  B states \cite{Hardy and Littlewood}
that the number of twins below a given bound
$N$ should be approximately equal to
\bee
\pi_d(N)\sim c_2 \int_2^N {du\over \ln^2(u)} = c_2{N\over \ln^2(N)}+\dots.
\label{H-L}
\ee
where the constant $c_2$ (sometimes called
``twin--prime'' constant)  %, see \cite{Wrench})
is defined in the following way:
\bee
c_2 \equiv 2 \prod_{p > 2} \biggl( 1 - {1 \over (p - 1)^2}\biggr) =
1.32032\ldots
\label{stalac2}
\ee

\begin{figure}[p]
\centering
\vspace{1.4cm}
\hspace{0.5cm}
\epsfxsize=13.5cm
\epsfbox{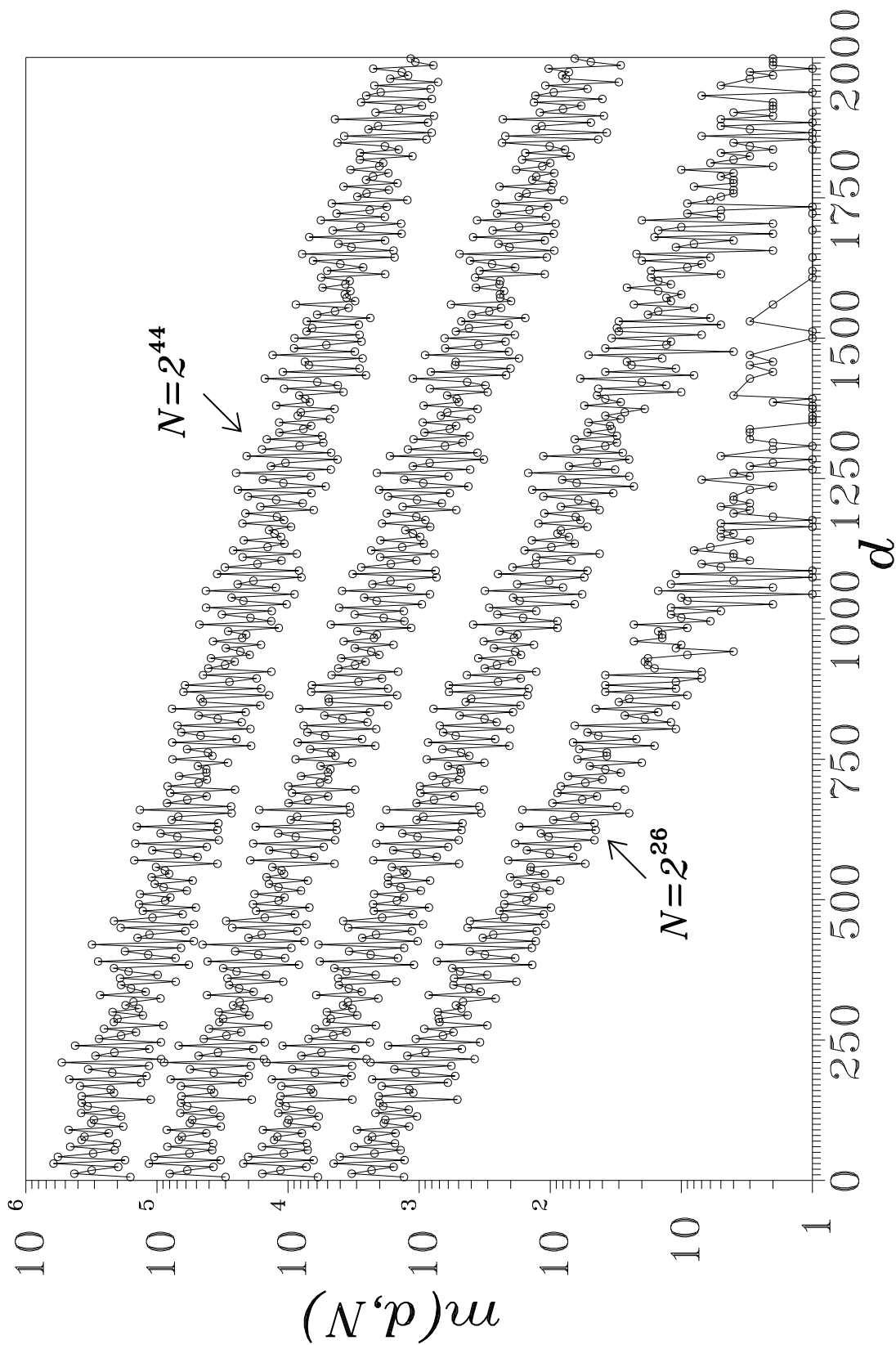}
\vspace{0.7cm}
\caption{The plot showing the dependence of histogram of gaps between
consecutive twins at $N=2^{26}, 2^{32}, 2^{38}$ and $2^{44}$. There is
linear scale on $x$ axis and logarithmical scale on $y$ scale.}
\label{twins}
\end{figure}

The problem of distribution of twins is very difficult. One of the
characteristics one can use for this purpose is the statistics of gaps
between consecutive twins. Let $d$ denote  the distance between two
consecutive twins measured as a arithmetical difference between last
primes constituting consecutive twins,  thus for example for twins
(29,31) and (17,19) $d=12$.  The distances can be only multiplicities of 6:
$d=6k$,  because all twins are of the  form $6k\pm 1$\cite{h}.
Next let $m(d,N)$ denote the number of
twins separated by $d$ and smaller than $N$:

$$ m(d,N)= {\rm number~ of~consecutive~ twins~}(p_m,p_{m+1}=p_m+2)~ {\rm and}
~(p_n,p_{n+1}=p_n+2)$$
\bee {\rm such~ that} ~ p_n-p_{m+1} = d, ~and~p_{n+1}< N
\label{mdN}
\ee
A few years ago I have obtained on the computer data for $m(d,N)$ for
$N=2^{26},\ldots N=2^{44}$ \cite{Wolf_tw}. A part of results is shown on the
Fig.\ref{twins}.

Recently there appeared two papers by P.F.Kelly and T.Pilling
\cite{Kelly1}, \cite{Kelly2}. They  have looked on the histogram of
separations between
consecutive twins measured by the {\it number of primes} in between, not just
the arithmetical difference as above.  Here I am providing
heuristic formulas which describe findings
of Kelly  and Pilling as well as much larger than in  \cite{Kelly1} and
\cite{Kelly2}
computer data which corroborates the obtained analytical relations.

\pagebreak

Let $\mu(s,N)$ denote the number of consecutive twins smaller than $N$
and separated by $s$ primes, with exception of primes constituting twins
in question, i.e. for example for  (5, 7) and (11, 13) $s=0$ while for
(17, 19) and (29, 31) $s=1$ (prime 23 lies in between). More strictly we
adopt the following definition:
%$$ \mu(s, N)~=~number~of~ consecutive~ twins~(p_m,p_{m+1}=p_m+2)~ and ~(p_m^\prime,
%p_{m+1}^\prime)$$
%\bee such~ that ~ \pi(p_m^\prime)-\pi(p_{m+1})-1 = s, ~and~p_{m+1}^\prime < N
%\ee
$$ \mu(s, N)~=~number~of~ consecutive~ twins~(p_m,p_{m+1}=p_m+2)~ and ~(p_n,
p_{n+1}=p_n+2)$$
\bee such~ that ~ \pi(p_n)-\pi(p_{m+1})-1 = s, ~and~p_{n+1}< N
\ee

Here $\pi(N)$ is as usual the number of primes$<N$ and
\bee
\pi(N)=\int_2^N {du \over \ln(u)} = {N\over \ln(N)}+\dots.
\label{pi_N}
\ee

The authors of \cite{Kelly1} have found that $\mu(s,N)$ for a given $N$
decreases exponentially with $s$. In fact these authors are working with
relative frequencies $\mu(s,N)/\sum_t \mu(t,N)$  but here I will use
absolute values  $\mu(s,N)$.

\begin{figure}[p]
\centering
\vspace{1.5cm}
\hspace{0.5cm}
\epsfxsize=13cm
\epsfbox{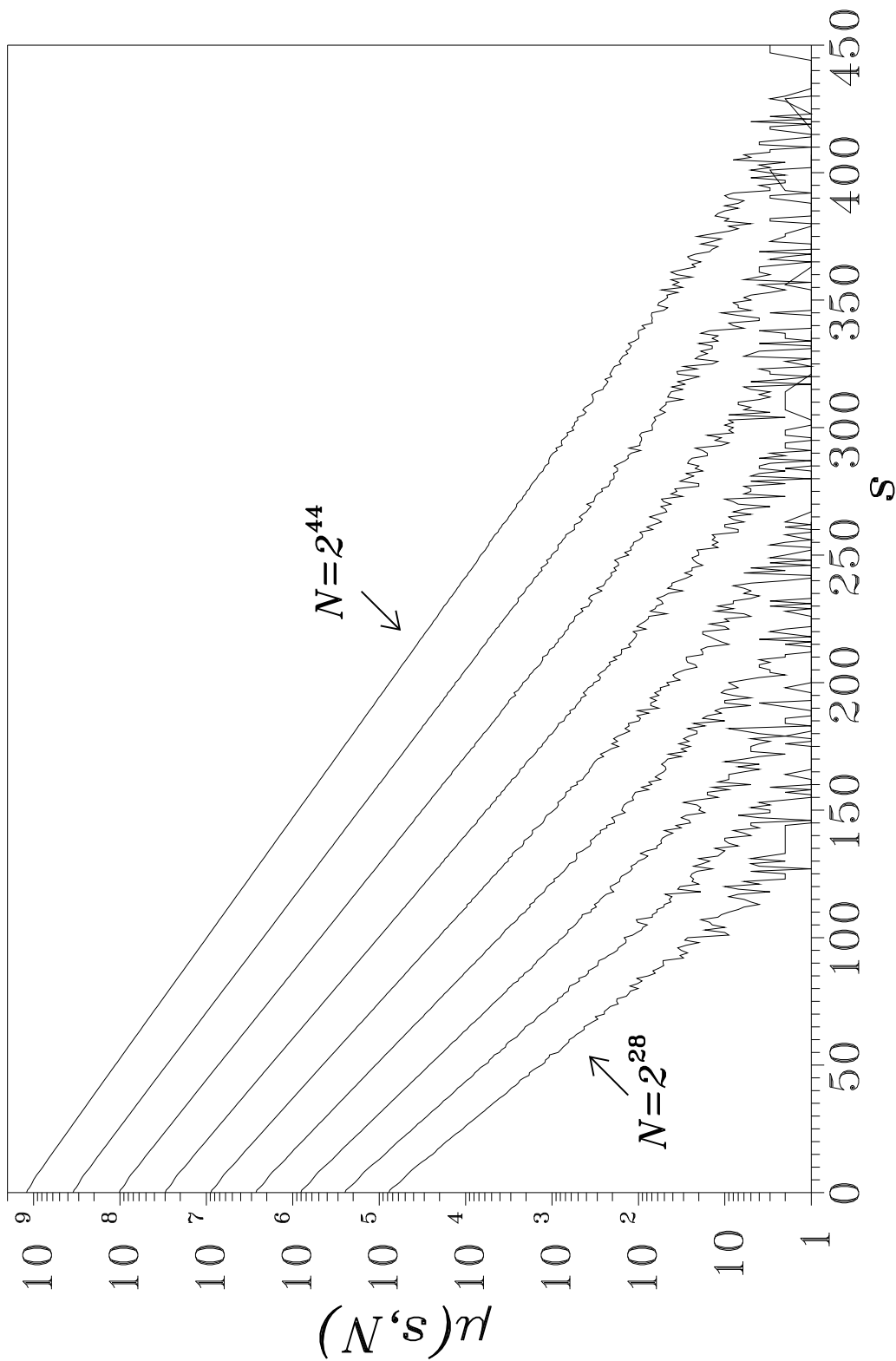}
\vspace{0.7cm}
\caption{The plot showing the dependence of the histogram  $\mu(s,N)$ on
$s$ at $N=2^{28}, 2^{30},
\ldots, 2^{44}$. There is a logarithmical scale on the $y$-axis , while on the
$x$-axis there is a linear scale.}
\label{mu}
\end{figure}

I have made the computer search up to $N=2^{44}\approx
1.76\times 10^{13}$ and counted the number of primes  between consecutive
twins. During the computer search the data representing the
function $\mu(s,N)$ were stored at values of $N$ forming the geometrical
progression with the ratio 4, i.e. at $2^{22}, \ldots, 2^{42},
2^{44}$. The resulting curves are plotted in the Fig.\ref{mu}.

Because the points lie on the straight lines on the semi-logarithmic
scale, we can infer from the Fig.\ref{mu} the ansatz\\
\bee
\mu(s,N) \sim A(N) e^{-B(N) s}.
\label{exponent}
\ee

The functions $A(N)$ and $B(N)$, giving the
intercepts and the slopes  of straight lines seen in the Fig.\ref{mu},
can be determined by exploiting two  identities %selfconsistency conditions
that $\mu(s,N)$ have to obey. First of all, the sum of $\mu(s,N)$ over all
$s$ is the total number of twins $<N$:
\bee
\sum_{s=0}^{s_{max}(N)} \mu(s,N)~=~\pi_2(N).
\label{self1}
\ee
Here $s_{max}(N)$ is the largest separation $s$ between twins $<N$, see
later discussion. The second selfconsistency condition comes from the
observation, that
\bee
\sum_{s=0}^{s_{max}(N)} s\mu(s,N) = \pi(N)-2\pi_2(N).
\label{self2}
\ee
In fact the above sum starts with $s=1$ what is important for the
relations obtained below. Putting the ansatz (\ref{exponent})  into
(\ref{self1}) and (\ref{self2}) and
collecting appropriately terms we end up with the geometrical series
with quotient $e^{-B(N)}$ and (\ref{self2})
is a differentiated

\pagebreak

\noindent geometrical series. Using the approximation
and summing in (\ref{self1}) and (\ref{self2}) up to infinity instead of
particular $s_{max}$ (which depends on $N$) I end up with the equations:
\bee
{A(N)\over 1-\exp\{-B(N)\}} = \pi_2(N)
\label{rozw1}
\ee
and
\bee
{A(N)e^{-B(N)} \over (1-\exp\{-B(N)\})^2} = \pi(N)-2\pi_2(N).
\label{rozw2}
\ee
From the Fig.\ref{mu} it is seen that the slopes $B(N)$ are smaller than 1 and
they decrease with $N$, e.g. $B(2^{42})=0.0503...,~B(2^{44})=0.0483...$.
Thus in the above relations I can put
$\exp\{-B(N)\}\approx 1 - B(N),~~B(N)\ll 1$  and solve the equations.
Finally the main conjecture is:
\bee
A(N)\approx{\pi_2^2(N) \over \pi(N)-2\pi_2(N)},
~~~~~~B(N)\approx{\pi_2(N) \over \pi(N)-2\pi_2(N)}.
\label{main_conj}
\ee

Putting here asymptotic forms of $\pi(N)$ and $\pi_2(N)$ I make the
following guess:
\bee A(N)\sim {c_2^2N^2 \over \ln^3(N)},~~~~~B(N)\sim {c_2 \over \ln(N)}.
\label{asympt}
\ee

The comparison of these formulae is given in Table 1, where the values
$A_{theor}(N)$  and $B_{theor}(N)$ obtained from (\ref{rozw1}) and
(\ref{rozw2}) are divided by values $A_{exp}(N)$  and $B_{exp}(N)$ obtained
by least-square method from data presented on Fig.\ref{mu}. The actual numbers
$A_{exp}$ and $B_{exp}$
depend on the number of points taken for linear regression method,
results in Table 1 are obtained when I skipped 15 first points and
40\% of last point, where large fluctuation appear
(in fact only 2  points determine straight line).
For the calculation of $A_{theor}(N)$  and $B_{theor}(N)$ I have used
{\it exact} values of $\pi(N)$ and $\pi_2(N)$ at $N=2^{22},\ldots
N=2^{44}$ from my earlier computer run \cite{Wolf}.

All ratios tend to 1 with increasing $N$, as it should be and asymptotic
values (\ref{asympt}) are smaller than exact values, since in (\ref{H-L})
and (\ref{pi_N}) dots denote positive terms which were skipped.

Another characteristic which can be used to test the conjecture
(\ref{main_conj}) is the question what are the largest separations
$s_{max}(N)$
between twins up to a given $N$. It corresponds to problem of the maximal gaps
between consecutive primes, what has a long history, see \cite{gaps}.
Heuristic argument is that the maximal gap appears only once and  hence it can
be obtained from the equation:

\bee
\mu(s_{max}(N),N)=1
\ee

\vskip 0.4cm
\begin{center}
{\sf TABLE {\bf I}}\\
\vskip 0.4cm
\begin{tabular}{|c|c||c||c||c||} \hline
$ N $ & $ A_{exp}(N)/A_{theor}(N) $ & $ B_{exp}(N)/B_{theor}(N) $  &
$ A_{exp}(N)/A_{asympt}(N) $ & $ B_{exp}(N)/B_{asympt}(N) $ \\ \hline$2^{22}$ &     0.964932 &     0.971223 &    1.528465 &     1.308701\\ \hline
$2^{24}$ &     0.983914 &     0.986198 &    1.475196 &     1.283024\\ \hline
$2^{26}$ &     0.967421 &     0.980584 &    1.389283 &     1.241010\\ \hline
$2^{28}$ &     0.972147 &     0.982905 &    1.353625 &     1.219132\\ \hline
$2^{30}$ &     0.965933 &     0.980277 &    1.312914 &     1.196643\\ \hline
$2^{32}$ &     0.970505 &     0.983126 &    1.291164 &     1.183385\\ \hline
$2^{34}$ &     0.973114 &     0.984629 &    1.270698 &     1.170806\\ \hline
$2^{36}$ &     0.975160 &     0.985884 &    1.252173 &     1.159566\\ \hline
$2^{38}$ &     0.976340 &     0.986722 &    1.235622 &     1.149610\\ \hline
$2^{40}$ &     0.977286 &     0.987373 &    1.220839 &     1.140667\\ \hline
$2^{42}$ &     0.977457 &     0.987681 &    1.206974 &     1.132435\\ \hline
$2^{44}$ &     0.978554 &     0.988279 &    1.195782 &     1.125445\\ \hline
\end{tabular} \\
\end{center}
\vskip 0.4cm

From the main conjecture (\ref{main_conj}) it follows that for large $N$:
\bee
s_{max}(N)\approx {\pi(N)\over \pi_2(N)}\biggl(2\ln\pi_2(N)
-\ln(\pi(N)-2c_2\pi_2(N))\biggr)
\label{s_max}
\ee

For large $N$ it goes into the
\bee
s_{max}(N)\sim {1\over c_2} \ln^2(N)
\ee
what differs from the Cramer's conjecture for maximal gaps between
consecutive primes by the factor $1/c_2$. It differs also by on power of
$\ln(N)$ from the asymptotic behavior of maximal gaps between
consecutive twins $G_2(N)\sim \ln^3(N)$
obtained just arithmetically: i.e. distance between twins
$p_m,p_{m+1}=p_m+2$ and $p_n, p_{n+1}=p_n+2$ is simply $p_n-p_m$ and not
$\pi(p_n)-\pi(p_{m+1}) -1$.
The comparison of actual values of $s_{max}(N)$ obtained from the
computer search are compared with formula (\ref{s_max}) in the Fig.3.

Finally I will make a remark on the problem of champions, i.e. the most
often occurring gaps. For prime numbers it was treated in
\cite{champion}. In \cite{Kelly2} the authors have made the assertion
that the champions for twins with gaps measured in terms of the primes
lying in between is always $s=0$. From the Fig.1 it is seen that when
gaps between twins are arithmetical differences then there will be a set
of champions consisting of $d=30$, next  emerging peak is $d=210$  and
so on, as it can be seen on the data presented in Fig.1. It will be
discussed in more detail in the forthcoming paper.

Finally let us notice that from the Fig. 2 it follows that the
separations between twins measured by the number of primes in between
follow exactly the Poissonian  behavior, see e.g. \cite{Poisson}.
Interestingly the change in the ``measuring sticks'' removes
oscillations from Fig.1 and leaves pure exponential decrease.

%{\bf Acknowledgment:} The author does not receive any support from the
%KBN or another source.

\vfill
\eject

\begin{figure}[p]
\centering
\vspace{1.4cm}
\hspace{0.5cm}
\epsfxsize=13cm
\epsfbox{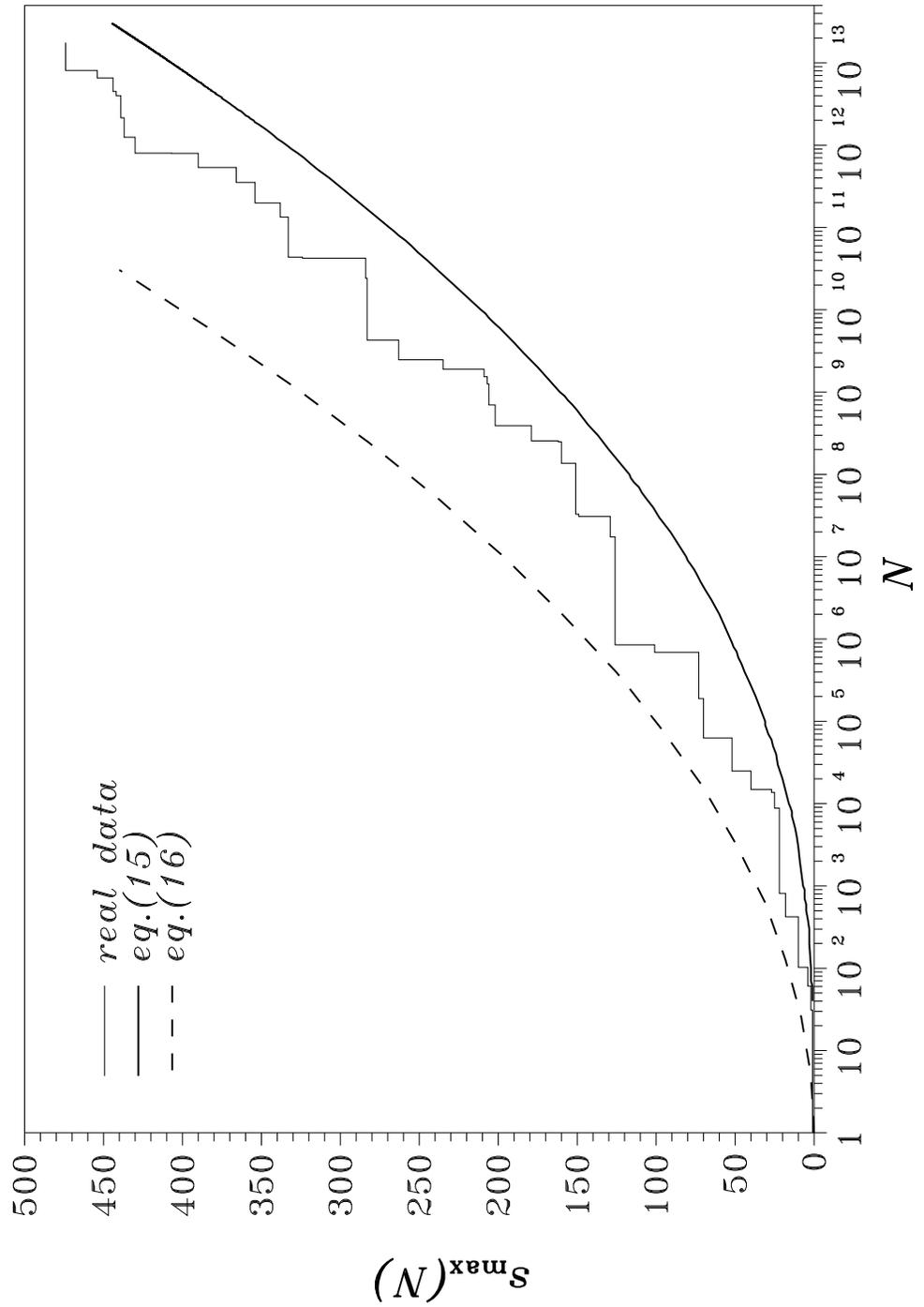}
\vspace{0.7cm}
\caption{The plot showing the dependence of $s_{max}(N)$ and comparison
with formula (\ref{s_max}).}
\label{fig_2}
\end{figure}
\pagebreak


\begin{thebibliography}{99}

\bibitem{Billingsley} P.Billingsley, ``Prime Numbers and Brownian
Motion'', {\it Amer. Math. Monthly} {\bf 80} (1973), pp.1099-1115

\bibitem{Hardy and Littlewood} G.H.Hardy and J.E. Littlewood,
{\it Acta Mathematica} {\bf 44} (1922), p.1-70

\bibitem{Wolf_tw} M.Wolf, 1999, unpublished

\bibitem{h} Cuesta-Dutari,-Norberto {\it Arithmetic of the sequences
$6n-1,\;6n+1$ and of twin primes.} Collect.-Math.
[Consejo-Superior-de-Investigaciones-Cientificas.-
Universidad-de-Barcelona.-Collectanea-Mathematica.-Seminario-Matematico-de-
Barcelona] 37 (1986), no. 3, p.211--227


\bibitem{Kelly1} P.F. Kelly and T.Pilling, {\it Characterization of the
Distribution of twin Primes},  arXiv:math.NT/0103191

\bibitem{Kelly2} P.F. Kelly and T.Pilling, {\it Implication of a New
Characterization of the Distribution of twin Primes}, arXiv:math.NT/0104205

\bibitem{gaps} H.Cramer, "On the order of magnitude of difference
between consecutive prime numbers", {\it Acta Arith.} {\bf 2} (1937),
p.23-46; D.Shanks, "On  Maximal Gaps between Succesive Primes",
{\it Math.Comp.} {\bf 18} (1964), p.464; J.H.Caldwell, ``Large Intervals
Between Consecutive Primes'',  {\it Math.Comp.}{\bf 25} (1971), p.909;
L.J.Lander and T.R.Parkin, ``On First Apperance of Prime
Differences'',  {\it Math.Comp.} {\bf 21} (1967), p.483; R.P.Brent,
``The First Occurrence of Certain Large  Prime Gaps'', {\it
Math.Comp.}{\bf 35} (1980), p.1435-1436;  J.Young and A.Potler, ``First
occurence Prime Gaps'', {\it Math.Comp.} {\bf 52} (1989), p.221-224;

\bibitem{Wolf} M.Wolf, {\it Some conjectures on the gaps between
consecutive primes},
preprint IFTUWr 894/95, available at http://www.ift.uni.wroc.pl/$\sim$mwolf


\bibitem{champion} A. Odlyzke, M. Rubinstein and M. Wolf, {\it Jumping
Champions}, Exp. Math., {\bf 8} (1999) 107-118.
Also available at http://www.ift.uni.wroc.pl/$\sim$mwolf and
http://www.research.att.com/$\sim$amo/doc/numbertheory.html

\bibitem{Poisson} O. Bohigas and M.J. Giannoni, in {\it Mathematical and
Computational Methods in Nuclear Physics}, edts. J.S. Dehesa et. al.
(Springer, Heidelberg, New York, 1984)


\end{thebibliography}
\end{document}